\documentclass[12pt]{amsart} 
\usepackage{amssymb,latexsym,epic}
\newlength{\fw}
\newlength{\jmr}
\newlength{\bernd}
\newlength{\ioan}
\newlength{\agk}
\newlength{\gkz}
\newtheorem{cor}{Corollary}
\newtheorem{lemma}{Lemma}

\newtheorem{van}{Vanishing Theorem for Resultants} 

\newtheorem{thm}{Theorem}
\newtheorem{rem}{Remark}	
\newtheorem{ex}{Example}
\newtheorem{dfn}{Definition}

\renewcommand{\qed}{$\blacksquare$} 
\newcommand{\relint}{\mathrm{RelInt}} 
\newcommand{\eps}{\varepsilon}
\newcommand{\cA}{\mathcal{A}}
\newcommand{\cD}{\mathcal{D}}
\newcommand{\cE}{\mathcal{E}}

\newcommand{\cH}{\mathcal{H}}

\newcommand{\cN}{\mathcal{N}}

\newcommand{\cP}{\mathcal{P}}

\newcommand{\conv}{\mathrm{Conv}}

\newcommand{\init}{\mathrm{in}}

\newcommand{\thth}{{\underline{\mathrm{th}}}}

\newcommand{\nd}{{\underline{\mathrm{nd}}}}
\newcommand{\Pro}{\mathbb{P}}
\newcommand{\Q}{\mathbb{Q}}
\newcommand{\R}{\mathbb{R}}
\newcommand{\C}{\mathbb{C}}
\newcommand{\N}{\mathbb{N}}
\newcommand{\Z}{\mathbb{Z}}
\newcommand{\bp}{\mathbf{\pi}}

\newcommand{\area}{\mathrm{Area}}

\newcommand{\divisor}{\mathrm{Div}}
\newcommand{\codim}{\mathrm{codim}}

\newcommand{\fan}{\mathrm{Fan}}

\newcommand{\res}{\mathrm{Res}}

\newcommand{\hilb}{\mathrm{\mathbf{Hilb}}}

\newcommand{\Zn}{\Z^n}
\newcommand{\Zno}{\Z^n\!\setminus\!\{\bO\}} 
\newcommand{\Qn}{\Q^n} 
\newcommand{\Rn}{\R^n}

\newcommand{\Cs}{\C^*}
 
\newcommand{\Ks}{K^*}
\newcommand{\Kn}{K^n}

\newcommand{\Qns}{\Qn\!\setminus\!\{\bO\}}

\newcommand{\Rns}{\Rn\!\setminus\!\{\bO\}} 

\newcommand{\Ksn}{{(K^*)}^n}

\newcommand{\cM}{\mathcal{M}}

\newcommand{\cC}{\mathcal{C}}
\newcommand{\cF}{\mathcal{F}}

\newcommand{\cT}{\mathcal{T}}
\newcommand{\bO}{\mathbf{O}}

\begin{document}
\title[Some New Applications of Toric Geometry]{Some New Applications of Toric 
Geometry} 

\author{J. Maurice Rojas}
\thanks{This research was partially funded by N.S.F. grant DMS-9022140 
and an N.S.F. Mathematical Sciences Postdoctoral Fellowship.}

\address{\hskip-\parindent
Massachusetts Institute of Technology\\
	 Mathematics Department\\
         77 Mass. Ave.\\
         Cambridge, MA \ 02139, U.S.A. } 

\date{\today}

\thanks{Revision of an earlier version which 
appeared in Foundations of Computational Mathematics, selected 
papers of a conference, held at IMPA in Rio de Janeiro, January 1997, 
Felipe Cucker and Mike Shub (eds.), Springer-Verlag. }

\begin{abstract}
This paper reexamines univariate reduction from a toric geometric 
point of view. We begin by constructing a binomial variant of 
the $u$-resultant and then retailor the generalized characteristic 
polynomial to fully exploit sparsity in the monomial structure 
of any given polynomial system. We thus obtain a fast new algorithm 
for univariate reduction and a better understanding of the underlying 
projections. As a corollary, we show that a refinement 
of Hilbert's Tenth Problem is decidable within single-exponential time. 
We also show how certain multisymmetric functions of the roots 
of polynomial systems can be calculated with sparse resultants. 
\end{abstract} 

\maketitle

\section{Introduction}
\label{sec:intro} 
We give a new approach for combatting degeneracy problems  
which occur when reducing large polynomial systems to univariate polynomials. 
Taking a positive attitude, we will actually make use of these degeneracies to 
better understand polynomial system solving. We do this by applying new 
observations involving toric varieties to establish faster, more 
reliable algorithms for univariate reduction.  Our techniques provide an 
intrinsic geometric setting for such reductions and fully exploit the sparsity 
of any polynomial system specified by its monomial term structure. 
 
Algebraically reducing problems involving polynomial systems  
to simpler questions involving a single univariate polynomial is a 
technique which has been known for over a century, and has been applied with 
increasing efficiency in computational algebra over the last two decades. 
(For example, 
see \cite{lazard,kaplak,pardokrick,fast} and the references therein.) 
However, there are difficulties with this method which {\em still}\/  
have not been completely addressed. Three such problems are the following: 
\begin{itemize}
\item[A:]{Since univariate reduction corresponds to 
projecting the solution space onto a projective line, 
the fiber over a point with finite coordinate might 
contain roots at infinity.} 
\item[B:]{The univariate polynomial one reduces to 
might be identically zero, due to (1) degeneracies in the underlying 
projection, or (2) the presence of infinitely many roots.} 
\item[C:]{Univariate reduction is still too slow for 
many problems of interest.} 
\end{itemize}
We circumvent A by the technique of {\em toric laminations}. By constructing 
a generalization of Canny's {\em generalized characteristic polynomial}\/ 
(GCP) \cite{gcp,mikebezout} --- the {\em toric}\/ (or {\em sparse}) GCP --- we 
can dispose of B$_2$. Our toric geometric approach also 
helps us better understand B$_1$. These two new techniques 
are contained in theorems 1--3 below. 

Although problem B$_2$ was solved earlier in \cite{renegar,
gcp}, the toric GCP greatly improves the complexity bounds given 
there. In particular, we address problem C by using the {\em 
sparse resultant}\/ \cite{gkz90,isawres} throughout our development. Fast 
methods for 
polynomial system solving were derived via the {\em sparse $u$-resultant}\/ in 
\cite{emiphd,isawres} but problems A and B were not addressed there.  
Our approach to univariate reduction thus unifies the GCP with the the sparse 
$u$-resultant. Furthermore, we generalize both techniques to certain {\em 
toric varieties}\/ \cite{kkms,tfulton,grobandpoly}.  

One advantage of a toric variety setting is reducing the problem of 
{\em exact}\/ root counting for polynomial systems, almost 
completely, to convex geometry. Generically sharp {\em upper bounds}\/ on the 
number of (complex) isolated roots, in terms of mixed volume, first appeared 
just over two decades ago \cite{kus75}. These bounds have since been 
generalized to algebraically closed fields \cite{dannie,toricint} and to 
various subsets of affine space (other than the algebraic torus) 
\cite{kho78,hsaff,toricint}. However, the following algorithm appears to be the 
first example of a convex geometric approach to counting the exact number 
of roots. (We will assume henceforth that all our polynomials and roots are 
considered over an algebraically closed ground field $K$.) 
\begin{thm} 
Let $F(x_1,\ldots,x_n)$ be an $n\!\times\!n$ polynomial system with support 
$E\!:=\!(E_1,\ldots,E_n)$ such that $\cM(E)\!>\!0$. 
Also let $P_E$ be the sum of 
the convex hulls of the $E_i$, and pick $a\!\in\!\Zno$ such that the segment 
$[\bO,a]$ is not parallel to any facet of $P_E$. Define 
$\bp_a(u_+,u_-)\!:=\!\res_{E_a}(F,u_+ +u_-x^a)$ where $\{u_\pm\}$ is a pair 
of algebraically independent indeterminates and $E_a\!:=\!(E,\{\bO,a\})$. 
Finally, let $\eps_\pm$ be the lowest exponent of $u_\pm$ occurring in any 
monomial of $\bp_a$. Then $F$ has {\em exactly} $\cN\!:=\!\cM(E)-\eps_+ 
-\eps_-$ isolated roots (counting multiplicities) in $\Ksn$, provided 
$\cN\!<\!\infty$. 
\end{thm} 
\noindent 
In the above, $\res_{\star}(\cdot)$ and $\cM(\cdot)$ respectively denote 
the sparse resultant and mixed volume \cite{isawres,mvcomplex,toricint}. 
The key contributions of the above theorem are (a) fibering the 
solution space by maps more general than coordinate projections, and (b) 
making explicit use of roots at {\em toric}\/ infinity when doing univariate 
reduction. 

We may also calculate the number of {\em distinct}\/ roots in a similar way. 
\begin{cor}
Following the notation and assumptions of theorem 1, let $\bp'_a(u_+,u_-)$ be 
the square-free part of $\bp_a(u_+,u_-)$. Also let $\eps'_\pm$ be 
the lowest exponent of $u_\pm$ occuring in $\bp'_a$. Then $F$ has 
{\em exactly} $\cN'\!:=\!\deg(\bp'_a)-\eps'_+-\eps'_-$ roots in $\Ksn$, 
provided $\cN'\!<\!\infty$ {\em and} the projection of roots 
$\{\zeta\!\in\!\Ksn \; | \; F(\zeta)\!=\!0\}\longrightarrow\{\zeta^a\}$ is 
injective. 
\end{cor} 
\noindent 
This corollary is proved alongside theorem 1 in section \ref{sec:chowtor}.  

By instead calculating a particular {\em coefficient}\/ of 
$\bp_a$, we can actually compute certain multisymmetric functions 
of the roots of $F$ and thus do more than just count the number of isolated 
roots. Admittedly, such a computation can be rather difficult, but recent 
advances in interpolation techniques, e.g., \cite{zip}, are making this 
approach increasingly feasible. So let us at least give explicit formulae 
for the coefficients of $\bp_a$ in terms of multisymmetric functions. 
\begin{dfn}
For any $w\!\in\!\Rns$, let $E^w_i$ be the set of points of $E_i$ 
having {\em minimal} inner product with $w$. Also, let 
$E^w\!:=\!(E^w_1,\ldots,E^w_n)$ and, when $w\!\in\!\Qns$, let 
$p_w\!\in\!\Zn$ be the first lattice point not equal to the origin encountered 
along the ray generated by $w$. Finally, let $L(E,a)$, the {\em 
ambiguity locus}, be the subvariety 
of $\cT_{P_E}\!\setminus\!\Ksn$ corresponding to the union of all faces of 
$P_E$ having an inner normal perpendicular to the segment $[\bO,a]$. 
\end{dfn}
\begin{rem}
The toric variety (over $K$), $\cT_P$, corresponding to a polytope $P$
is detailed in \cite{ksz92,tfulton,gkz94,grobandpoly,toricint}. Polynomial
roots in a toric compactification are described (in a manner closest
to our present framework) in \cite{toricint}. 
\end{rem}
\begin{cor}
\label{cor:multi}
Following the notation of theorem 1, assume instead that $a$ is {\em
any} point in $\Zno$ and $\eps_+,\eps_-\!<\!\infty$. 
Then $\bp_a$ is a homogeneous polynomial of degree $\cM(E)$ and, for 
any $d\!\in\!\{0,\ldots,\cM(E)-\eps_+-\eps_-\}$, the 
coefficient of $u^{d+\eps_-}_-u^{\cM(E)-d-\eps_+}_+$ in $\bp_a$ is precisely 
\[C \sum \prod\limits^d_{j=1} \zeta(j)^a, \]
where $C\!\in\!\Ks$ is independent of $d$, the sum ranges over all cardinality 
$d$ multisets $\{\zeta(1),\ldots,\zeta(d)\}\!\subseteq\!Z_F$, 
and $Z_F$ is the multiset of roots of $F$ in $\Ksn\cup L(E,a)$. Furthermore, 
if $\eps_+\!=\!\eps_-\!=\!0$ then $C$ is precisely 
\[ (-1)^n\left(\prod\limits_w \res_{E^w}(F)^{-p_w\cdot a}\right), \] 
where the product ranges over all inner facet directions of $P_E$ having 
negative inner product with $a$. \qed  
\end{cor} 
\noindent 
Multisets are simply sets where repeated elements are allowed, and whose 
subsets have repetitions 
appropriately restricted \cite{gkp}. For example, if $Z_F$ contains $-1$ as a 
root of multiplicity $5$, then no subset of $Z_F$ can have more than five 
$-1$'s. So in essence, the coefficients of $\bp_a$ (when suitably normalized) 
are just elementary symmetric functions of $a$-monomials of certain roots of 
$F$. 

The above corollary follows easily from theorem 1 and the 
Pedersen-Sturmfels product formula \cite{chowprod}.\footnote{We 
also need the fact that the Pedersen-Sturmfels formula, originally 
stated only over $\C$, remains true over a general algebraically 
closed field. This is proved in \cite{resvan}.} Moreover, since 
the sparse resultant is invariant under muliplication of the individual 
polynomials by a monomial, we can compare the coefficients of 
$\bp_a$ and $\bp_{-a}$ to obtain the following beautiful identity.\footnote{
Alicia Dickenstein pointed out the $K\!=\!\C$ of the above formula (due to 
Cattani, Dickenstein and Sturmfels) to the author 
at the conference where this paper was presented. Shortly after, 
the author realized that this formula could also follow easily from 
corollary \ref{cor:multi}.}
\begin{cor}
Following the notation and assumptions of corollary \ref{cor:multi}, 
$\eps_+\!=\!\eps_-\!=\!0 \Longrightarrow$
\[ \prod\limits_{\zeta\in Z_F} (\zeta^a)= \prod\limits_w 
\res_{E^w}(F)^{p_w\cdot a} \]
where the product ranges over {\em all} inner facet directions of $P_E$. 
\qed
\end{cor}

Recall that a {\em ridge}\/ is a polytope face of codimension $2$ 
\cite{ziegler}. Although the ``bad'' case of theorem 1 ($\cN\!=\!\infty$) 
is an annoyance, it is a small annoyance and can actually be made use of in 
certain cases. 
\begin{thm}
\label{thm:ambi}
Following the notation and assumptions of theorem 1, 
$\cN\!=\!\infty \Longleftrightarrow F$ has a root lying in $L(E,a)$ {\em 
or} infinitely many roots in $\Ksn$. Also, within the space 
of all systems with support contained in $E$, the $F$ with $\cN\!=\!\infty$ 
form a subvariety of codimension $\geq\!2$. 
\end{thm}
\noindent
Thus, if the roots of $F$ avoid the ambiguity locus, $F$ has infinitely 
many roots in $\Ksn \Longleftrightarrow \cN\!=\!\infty$.
In particular, it follows directly from the definition of $L(E,a)$ that 
$\codim L(E,a)\!=\!2$, under the assumptions of theorem \ref{thm:ambi}. 
(Since $L(E,a)$ then corresponds to a finite set of ridges \cite{tfulton,
toricint}.)

We can completely avoid amiguity loci by replacing the binomial 
$u_+ +u_-x^a$ in our preceding results with another specially chosen 
indeterminate polynomial. This is detailed further in \cite{tgcp} and, to a 
certain extent, in the result we quote below. However 
it is worth noting that the roots of $F$ {\em always}\/ avoid $L(E,a)$ when 
$n\!<\!3$ and $E$ is the support of $F$ (cf.\ lemma \ref{lemma:once}). 
\begin{rem}
The existence of infinitely many roots of $F$ at 
toric infinity and the existence of a root of $F$ within the ambiguity 
locus are independent in general (cf.\ section \ref{sec:ex}). Thus it is 
worthwhile to know something about $L(E,a)$.
\end{rem} 
\begin{rem}
\label{rem:egg}
A correspondence useful for visualizing theorems 1 and 2 is to 
identify toric infinity ($\cT_P\!\setminus\!\Ksn$) with the 
boundary of the polytope $P$ (cf.\ definition \ref{dfn:toric}). 
Computing $\bp_a$ then amounts to splitting $\partial P_E$ into two signed 
halves, much like cracking an 
egg-shell. In particular, the half with sign $\pm$ consists precisely of 
those faces of $P_E$ having an inner normal with $\pm$ inner product with 
$a$. The ambiguity locus is then precisely the common boundary of these 
two halves. We will see later that $\bp_a$ can be viewed as a coordinatization 
of a particular rational map $\cT_{P_E} \dashrightarrow \Pro^1_K$. Thus, 
provided no roots lie in the ambiguity locus, $\eps_{\pm}$ is the sum of the 
intersection numbers of the roots of $F$ which lie in the $\pm$ half.  
\end{rem}

A general way to deal with the existence of {\em 
infinitely}\/ many roots of $F$ within a given toric compactification 
is the following result announced below. 
\begin{thm} \cite{tgcp}  
Following the notation of theorem 1, suppose $g(x)\!:=\!\sum_{e\in A} 
u_e x^e$ where the $u_e$ are algebraically independent indeterminates and 
$A\!\subset\!\Zn$ is nonempty and finite. Suppose further that 
$\cP\!:=\!(P_1,\ldots,P_n)$ is an $n$-tuple of integral polytopes in $\Rn$ 
containing the support of $F$, and let $D$ be an irreducible fill of $\cP$. 
Setting $F^*\!:=\!(\sum_{e\in D_i} x^e \; | \; i\!\in\![1..n])$, define 
$\cH(s;u)\!:=\!\res_{(E,A)}(F-sF^*,g)$, where $s$ is a new indeterminate. 
Finally, considering $\cH$ as a polynomial in $s$ with 
coefficients in $K[u]$, let $\cF_A(u)$ be the coefficient of the lowest power 
term of $\cH$. Then $\cF_A$ is a homogeneous polynomial with the following 
properties:
\begin{enumerate}
\item{$\cF_A$ is divisible by $g(\zeta)$, for any root $\zeta\!\in\!\Ksn$ of 
$F$.} 
\item{If $\cP$ is compatible with $\conv(A)$, then $\cF_A$ has degree 
$\cM(\cP)$ and to every irreducible factor of $\cF_A$ there naturally 
corresponds a root of $F$ in $\cT_P$, where $P\!:=\!\sum P_i$.}
\end{enumerate}
We call $\cH$ a {\em toric} (or {\em sparse}) {\em generalized 
characteristic polynomial for $(F,A)$.} \qed  
\end{thm} 
The combinatorial portion of the above 
theorem can be explained as follows: We say that $D$ {\em fills}\/ $\cP$  iff 
$D\!:=\!(D_1,\ldots,D_n)$ satisfies $D_i\!\subseteq\!P_i$ for all 
$i\!\in\![1..n]$ and  $\cM(D)\!=\!\cM(\cP)$ \cite{convexapp,rojaswang}. An {\em 
irreducible}\/ fill is then simply a fill which is minimal with respect to 
$n$-tuple containment. Finding such a $D$ amounts to a combinatorial 
preprocessing step which need only be done once for a given set of 
problems, provided $E$ is fixed. Compatibility is defined in section 
\ref{sec:chowtor} and the toric GCP is briefly compared 
to the original GCP in section \ref{sec:toricgcp}.
Combined with the toric variety version of 
Bernshtein's theorem \cite{toricint}, we see that the above theorem gives 
us a resultant-based method to work with the isolated roots of a sparse 
polynomial system, even 
in the presence of infinitely many roots. Furthermore, the roots 
defined by $\cF_A$ which do {\em not}\/ correspond to any isolated root 
of $F$ can be usefully interpreted as contributions from the excess components. 

We apply our main theorems to some simple examples in section \ref{sec:ex}. 
Our results are also very useful for root 
counting over $\R$, $\Q$, and even $\Z$. For example, one can combine standard 
univariate techniques such as Sturm sequences \cite{glrr,marie1} with 
our preceding results to quickly count the number of real roots \cite{nash}. 
Better still, we can go to an even smaller ring: In section \ref{sec:hilbert} 
we will give a short proof (under some mild hypotheses) of the following 
refinement of Hilbert's Tenth Problem \cite{maty}. 
\begin{thm} \cite{resvan} 
\label{thm:hilbert}
Consider, for any nonnegative integer $d$, the  
following collection of Diophantine systems:\\
{${\bf \hilb(d)}${\bf :}} Multivariate polynomial systems with integer 
coefficients and $d$-dimensional {\bf complex} solution set.\\ 
Then there is an algorithm which, given any instance of $\hilb(0)$, 
finds all integer solutions, or certifies that there are none, within 
single-exponential time. \qed 
\end{thm} 

All of the above results stem from an observation on the 
vanishing of the sparse resultant. In particular, one basic 
property of the sparse resultant is that 
\[(\star) \ \ \ \ \ \ \ \ \ (F,f_{n+1}) \mathrm{\ has \ a \ root \ in \ } \Ksn 
\Longrightarrow \res_{\bar{E}}(F,f_{n+1})\!=\!0\] 
(provided $\bar{E}\!=\!(E,E_{n+1})$ and $f_{n+1}$ is a polynomial over $K$ with 
support contained in $E_{n+1}\!\subset\!\Zn$). However, the converse assertion 
is \underline{{\bf not}}\/ always true. (The polynomial system 
$(x+y+1,x+y+2,x+y+3)$ with $\bar{E}\!=\!\{\bO,\hat{e}_1,\hat{e}_2\}^3$ 
gives a simple counter-example reducing to the vanishing of a $3\times 3$ 
determinant.) 
The correct statement in which {\em both}\/  
implications hold seems to be known only folklorically, while the {\em 
unmixed}\/ case $E_1\!=\!\cdots\!=\!E_{n+1}\!=\!\cA$ (over $\C$) is contained 
in various recent works, e.g., \cite{ksz92,gkz94}. So we will make use of the 
following more general result. 
\begin{van} \cite{resvan} 
Suppose $f_i$ is a polynomial over $K$ with support {\em contained in}\/ 
$E_i\!\subset\!\Zn$ for all $i\!\in\![1..n+1]$. Then, provided 
\[\cM(E_1,\ldots,\widehat{E_i},\ldots,E_{n+1})>0\] 
for some $i\!\in\![1..n]$, 
\[\res_{\bar{E}}(f_1,\ldots,f_{n+1})\!=\!0 \Longleftrightarrow 
\bigcap\limits^{n+1}_{i=1} \cD_{P_{\bar{E}}}(f_i,\conv(E_i))\!\neq\!
\emptyset,\]
where $\bar{E}\!:=\!(E_1,\ldots,E_{n+1})$ and 
$P_{\bar{E}}\!:=\!\sum\limits^{n+1}_{i=1} \conv(E_i)$. \qed 
\end{van}
\noindent
In the above, $\cD_P(f_i,P_i)$ is the {\em toric effective 
divisor}\/ of $\cT_P$ corresponding to $(f_i,P_i)$ (cf.\ definition 
\ref{dfn:toric}). 
This result provides a geometric analogue, over a {\em general}\/ 
algebraically closed field, of the product formula for the sparse resultant 
\cite{chowprod}. 

We close this introduction with an important note: 

In addition to the applications we have presented, our approach to 
univariate reduction provides a first step toward 
an intrinsic method for root finding in the algebraic torus $\Ksn$. 
Although one could in principle reduce to various         
coordinate subspaces to perform ``sparse'' elimination theory in $\Kn$, such 
an approach 
runs into complications with intersection multiplicities on the coordinate 
hyperplanes and is really {\em not}\/ intrinsic to $\Kn$. We propose a more 
intrinsic (and efficient) approach to {\em affine elimination theory}\/ by 
defining a new resultant operator --- the {\em affine}\/ sparse
resultant --- in \cite{aff}. Following this route,
{\em all}\/ of our main theorems generalize easily to
root finding in affine space minus an arbitrary union of
coordinate hyperplanes. This is pursued in greater depth in
\cite{toricint,aff} and the complexity bound of
theorem \ref{thm:hilbert} is further refined in the latter work.

\section{Related Work}
We point out a few important recent approaches to elimination 
theory which can profit from our results. 

From an applied angle, our observations on degeneracies and handling polynomial 
systems with infinitely many roots nicely complement the work of Emiris and 
Canny \cite{isawres}. In particular, their sparse resultant based algorithms 
for polynomial system solving can now be made to work even when  
problem B occurs. Also, an added benefit of working torically (as 
opposed to the classical approach of working in projective space) is the 
increased efficiency of the sparse resultant: the resulting matrix calculations 
(for polynomial system solving) are much smaller and faster. In particular, 
whereas it was remarked in \cite{gcp} that Gr\"obner basis methods are likely 
to be faster than the GCP for sparse polynomial systems, the toric GCP 
appears to be far more competitive in such a comparison. 

From a more theoretical point of view, our focus on univariate reduction is a 
useful addition to Sturmfels' foundational work on sparse elimination theory 
in the large \cite{sparseelim}. (We also point out that \cite{sparseelim} 
provides a wonderfully clear introduction to some of the toric variety 
techniques we refer to.) 

Going further into the intersection of algebraic geometry and 
complexity theory, there is a beautiful new approach to elimination theory 
founded by the school of Heintz, et.\ al. \cite{pardokrick,fast}. Our results  
fill in some toric geometry missing from their more algebraic framework. 
In particular, we propose that toric laminations, i.e., 
``fibering'' $\cT_{P_E}$ by {\em toric}\/ projections (instead of fibering 
$\Ksn$ by coordinate projections), can enhance their methods. For example, 
although their formulation in terms of straight-line programs is extremely 
general, it 
appears that their complexity bounds can be significantly improved with our 
techniques in the case of polynomial systems which are ``sparse'' in the 
specified supports sense we use here. This is explored further in 
\cite{resvan}.  
 
\section{Examples}  
\label{sec:ex} 
We will give a $2\times 2$ example of our first two theorems. 
 
\subsection{Root Counting with Sparse Resultants} 
Let $n\!=\!2$ and consider the bivariate polynomial system 
$F\!:=\!(x^3+y^4-1,x^4+y^5-1)$. Also let $E_1$ and $E_2$ respectively be the 
supports $\{\bO,(3,0),(0,4)\}$ and $\{\bO,(4,0),(0,5)\}$. 
Clearly then the Newton polygons are two triangles and the segment 
$[\bO,(1,1)]$ is not parallel to any edge of the quadrilateral 
$P_E\!:=\!\conv(\bO,(7,0),(0,9),(4,4))$. So we can 
set $a\!:=\!(1,1)$ and apply theorem 1 to count the roots of $F$.

We are now ready to compute the sparse resultant $\res_{E_a}(F,u_++u_-xy)$, 
where $E_a\!:=\!(E_1,E_2,\{\bO,(1,1)\})$ and the 
coefficients of $F$ have been replaced by algebraically independent 
indeterminates.\footnote{The published version of this paper contains an 
erroneous computation (due to a software error) of this resultant. This 
was observed by John Dalbec at the conference where this paper was presented, 
and the calculation below is due to him. The author is deeply indebted to 
Professor Dalbec's astute observation.} To do this, we can use {\tt Maple} to 
calculate the Sylvester ($2\times 1$) resultants of the pairs 
$(f_2,u_++u_-xy)$ and 
$(f_1,u_++u_-xy)$, eliminating the variable $y$. Let us respectively call 
these resultants $f_{-1}$ and $f_{-2}$. Using the Sylvester resultant one last 
time on $(f_{-1},f_{-2})$ to eliminate $x$, we arrive at a {\em multiple}\/ 
of the sparse resultant we are looking for. (This fact follows 
easily from the characterization of the sparse resultant as the 
(nonzero) polynomial of lowest total degree satisfying property 
($\star$) from the introduction \cite{sparseelim,resvan}.) 

The preceding computations, modulo typing, takes only fractions of a second 
on a Sun 4m Computeserver. From \cite{combiresult} we know that the sparse 
resultant itself should have degree $\cM(E_1,E_2)+\cM(E_1,\{\bO,a\})+
\cM(E_2,\{\bO,a\})\!=\!16+7+9\!=\!32$. (This calculation is easily done 
by hand via the general formula $\cM(E_1,E_2)\!=\!\area(\conv(E_1+E_2))-
\area(\conv(E_1)) -\area(\conv(E_2))$ \cite{schneider,mvcomplex}.) So by 
factoring with {\tt Maple} (which takes just under $2$ seconds) we can isolate 
the sparse resultant and specialize coefficients to obtain 
\[ \bp_{(1,1)}(u_+,u_-) = \res_{E_a}(x^3+y^4-1,x^4+y^5-1,u_+ +u_- xy)  =  \]  
\[20u^9_- u^7_+ + 31u^8_- u^8_+ + 12u^7_- u^9_+ + 14u^5_- u^{11}_+ 
+14u^4_- u^{12}_+ + 7u^3_- u^{13}_+ -9u^2_- u^{14}_+ + u_- u^{15}_+ +  u^{16}_+  \]
\noindent 
So $\eps_+\!=\!7$, $\eps_-\!=0$, and theorem 1 tells us that $F$ 
has exactly $16-7\!=\!9$ roots, counting multiplicities, in $(\Cs)^2$. 
Furthermore, a one-line {\tt Maple} calculation shows that $\bp_{(1,1)}$ is 
irreducible over $\Q$. So corollary 1 (and basic Galois theory) tells us that 
the multiplicity of every root of $F$ in $(\Cs)^2$ is $1$. 

The roots of $F$ can also be counted via Gr\"obner basis methods but the toric 
geometry in this example is rather nice 
to observe. In particular, the coordinate cross $\{xy\!=\!0\}$ is naturally 
embedded in $\cT_{P_E}$ and corresponds  precisely to the two lower-left edges 
of $P_E$ (cf.\ remark \ref{rem:egg} and definition \ref{dfn:toric}). 
Furthermore, the ambiguity locus consists of 
two distinct points which, due to the fact that $\bp_a$ is not identically 
$0$ (cf.\ theorem 2), the roots of $F$ have luckily avoided. Also 
$\eps_+\!=\!7 \Longrightarrow 
F$ has precisely $7$ roots, counting multiplicities, on the coordinate cross. 
The latter fact can also be easily checked with a little commutative algebra.

We have included timing information solely for illustrative purposes. 
In particular, our calculations can be sped up tremendously 
with suitably specialized code. 


\section{A Refinement of Hilbert's Tenth Problem}
\label{sec:hilbert} 
We will present a concise proof of theorem 4 under the 
following three hypotheses:
\begin{enumerate}
\item{The instances of $\hilb(0)$ considered are 
$n\times n$ polynomial systems.}
\item{The complex roots of an instance must have all 
coordinates nonzero.}
\item{No instance can have a complex root at toric 
infinity.} 
\end{enumerate} 
Note that (1) and (2) are of a combinatorial nature: 
they restrict the number and dimensions of the supports, 
as well as their intersections with the coordinate subspaces. 
On the other hand, (3) is of a more algebraic nature and, 
when the supports are fixed, is easily shown to fail only on a 
codimension $1$ subvariety of polynomial systems, e.g., \cite[Main Theorem 
2]{toricint}. 

The above hypotheses are present mainly for technical reasons and 
are removed in \cite{resvan}. They may also be removed with 
other recent techniques \cite{bpr,pardokrick,fast}, but through 
our toric techniques it is possible to obtain a general complexity 
bound with near {\em quadratic}\/ dependence on the mixed volume
\cite{tgcp}. This would be the best asymptotic bound to date for solving 
Diophantine systems. 

\noindent
{\bf Proof of Theorem 4:} First we outline our algorithm: 
\begin{itemize}
\item[Step 1:]{Following the notation of theorem 3, let $A$ 
be the vertices of a standard $n$-simplex in $\Rn$ and $\cP$ 
the $n$-tuple of Newton polytopes of $F$. Form the polynomial $\cF_A$.} 
\item[Step 2:]{Using the Lenstra-Lenstra-Lovasz algorithm 
\cite{lll}, factor $\cF_A$ over $\Z$.} 
\item[Step 3:]{Any integral root $(z_1,\ldots,z_n)\!\in\!\Zn$ of $F$ 
corresponds precisely to an integral factor of $\cF_A$ of the form  
$u_\bO+z_1u_{\hat{e}_1}+\cdots+z_nu_{\hat{e}_n}$. } 
\end{itemize}

Step 1 runs in time polynomial in the mixed volume, applying hypotheses 
(1) and (3), and the results of \cite{sparseelim,cxelim}. In fact, we are 
merely computing a ``sparse'' variant of the $u$-resultant in Step 1. To prove 
this (and the fact that this sparse $u$-resultant is not identically 
$0$) one uses slightly more general versions of theorems 1 and 2 
\cite{resvan}. One also needs the fact that hypotheses (1) and (3) together 
imply that the complex zero set of $F$ in $\Ksn$ is zero-dimensional, but this 
is an immediate consequence of \cite[Corollary 3]{toricint}. We also point 
out that working more generally in terms of $\cF_A$ aids in later removing 
hypothesis (3) \cite{resvan}. 

To see that Step 2 runs within time single-exponential in the input size, one 
need only observe that (a) the coefficient growth in the sparse $u$-resultant 
calculation can be reasonably bounded, and (b) the LLL-algorithm runs within 
time single-exponential in $n$ and polynomial in the bit-sizes of the 
coefficients of $F$. Part (a) is detailed in \cite{resvan}, and part (b) 
follows easily from the univariate case covered in \cite{lll}. 

That we can terminate as described in Step 3 follows immediately 
from hypothesis (2) and the Vanishing Theorem for Resultants. 
\qed 

We conclude this section with some brief notes on related work.

First, we point out the rather surprising fact that restructuring 
Hilbert's Tenth Problem as the union $\bigcup^\infty_{d=0} \hilb(d)$ seems to 
be new. (In particular, Hilbert's Tenth Problem was originally 
stated as deciding the existence of a {\em single}\/ integral root for 
{\em one}\/ polynomial in several variables.) Although the results of  
\cite{jones} imply that $\hilb(8)$ is undecidable, nothing seems to 
be known about $\hilb(d)$ for $0<d<8$. We find this 
shocking, considering the vast effort within arithmetic geometry to use 
complex geometric invariants in Diophantine problems, e.g., Falting's theorem 
and the deep conjectures of Lang and Vojta \cite{arithgeom}. 

We also point out that $\hilb(0)$ being a single-exponential problem 
may not be so new: Teresa Krick and Luis-Miguel Pardo-Vasallo have 
pointed out to the author that theorem 4 can also be readily 
derived from \cite{pardokrick,fast}. (Indeed, such an argument 
could have also been derived directly from \cite{lll} and \cite{renegar} 
in 1987.) Pardo-Vasallo has also pointed out that Marie-Fran\c{c}oise 
Roy lectured in 1995 on a similar result, presumably derived 
from \cite{bpr}. Nevertheless, optimal bounds are open territory; not 
to mention the $6$-long gap in the status of $\hilb(d)$. 

\section{Toric Varieties and Sparse Resultants} 
\label{sec:chowtor} 
Our notation is a slight variation of that used in 
\cite{tfulton}, and is described at greater length in 
\cite{toricint}. 
We will assume the reader to be familiar with normal fans of polytopes and  
the construction of a toric variety from a fan \cite{tfulton,gkz94}.
However, we will at least list our cast 
of main characters:
\begin{dfn} \cite{toricint} 
\label{dfn:toric}
Given any $w\!\in\!\Rn$, we will use the following notation:
\begin{itemize}
\item[$T=$]{The algebraic torus $\Ksn$}
\item[$P^w=$]{The face of $P$ with inner normal $w$} 
\item[$\sigma_w=$]{The closure of the cone generated by the inner normals of 
$P^w$} 
\item[$U_w=$]{The affine chart of $\cT_P$ corresponding to the cone
$\sigma_w$ of $\fan(P)$} 
\item[$L_w=$]{The $\dim(P^w)$-dimensional subspace of $\Rn$ parallel to $P^w$} 
\item[$x_w=$]{The point in $U_w$ corresponding to the semigroup
homomorphism $\sigma^\vee_w\cap\Zn \longrightarrow \{0,1\}$ mapping $p
\mapsto \delta_{w\cdot p,0}$, where $\delta_{ij}$ denotes the Kronecker
delta}  
\item[$O_w=$]{The $T$-orbit of $x_w=$ The $T$-orbit corresponding to $\relint 
P^w$} 
\item[$\cE_P(Q)=$]{The $T$-invariant Weil divisor of $\cT_P$ corresponding to 
a polytope $Q$ with which $P$ is compatible} 
\item[$\divisor(f)=$]{The Weil divisor of $\cT_P$ defined by 
a rational function $f$ on $\Ksn$} 
\item[$\cD_P(f,Q)=$]{$\divisor(f)+\cE_P(Q)=$ The toric effective divisor of 
$\cT_P$ corresponding to $(f,Q)$} 
\item[$\cD_P(F,\cP)=$]{The (nonnegative) cycle in the Chow ring of 
$\cT_P$ defined by $\bigcap^n_{i=1} \cD_P(f_i,P_i)$ } 
\end{itemize}
\end{dfn}

We will say that a polytope $P$ (resp.\ an $n$-tuple $\cP$) is 
{\em compatible}\/ with $Q$ iff every cone of $\fan(Q)$ 
is a union of cones of $\fan(P)$ (resp.\ a union of cones of $\fan(P_i)$ 
for each $i$) \cite{khocompat,toricint}. 
We will also make frequent use of the natural correspondence 
between the face interiors $\{\relint P^w\}$ and the $T$-orbits 
$\{O_w\}$ \cite{ksz92,tfulton,gkz94}. The following lemma gives a more 
explicit algebraic analogy between the vertices of $P$ and the maximal affine 
charts of $\cT_P$.  
\begin{lemma} \cite{toricint}
\label{lemma:once}
Suppose $F$ is a $k\times n$ polynomial system over $K$ with support 
contained in a $k$-tuple of integral polytopes $\cP\!:=\!(P_1,\ldots,P_k)$ in 
$\Rn$. Assume further that $P$ is a rational polytope in $\Rn$ compatible with 
$\cP$. Then the defining ideal in $K[x^e \; | \; e\in\sigma^\vee_w\cap \Zn]$ of 
$U_w\!\cap\!\cD_P(F,\cP)$ is $\langle x^{b_1}f_1,\ldots,x^{b_k}f_k \rangle$, 
for any $b_1,\ldots,b_k\!\in\!\Zn$ such that $b_i+P^w_i\subseteq L_w$ for all
$i\!\in\![1..k]$. \qed
\end{lemma} 

Theorems 1 and 2 will follow easily from the Vanishing Theorem 
for Resultants via following two lemmata below.
\begin{lemma}
\label{lemma:chow}
Following the notation of theorem 1, define $$\cP_E\!:=\!(\conv(E_1),\ldots,
\conv(E_n)).$$ Then the map $x \mapsto x^a$ extends to a 
proper morphism $\phi_a : \cT_{P_E}\!\setminus\!L(E,a) \longrightarrow 
\Pro^1_K$. Also, if $\cD_{P_E}(F,\cP_E)$ has 
zero-dimensional intersection with $\Ksn$ and avoids the ambiguity locus, 
then $\bp_a$ is, up to a nonzero constant multiple, the Chow form of the 
subvariety $\phi_a(\cD_{P_E}(F,\cP_E))$ of $\Pro^1_K$. \qed 
\end{lemma} 
\begin{lemma}
\label{lemma:forbid}
Following the notation of theorem 1, pick points 
$a(1),\ldots,a(n\!-\!1)\!\in\!\Zn$ 
which together generate a hyperplane in $\Rn$ which is not parallel 
to any facet of $P_E$. Then $F$ has infinitely many roots in $\Ksn 
\Longrightarrow$ the polynomial $p(u_+,u_-)\!:=\!\prod^{n-1}_{i=1} \bp_{a(i)}$ 
is identically $0$. \qed   
\end{lemma}
\noindent 
Sufficiently armed, let us now prove theorem 2. 

\noindent 
{\bf Proof of Theorem 2:} Suppose that $F$ has no roots within $L(E,a)$ and 
only finitely many roots within $\Ksn$. Then it follows from 
lemma \ref{lemma:chow} that $u_+$ 
(resp.\ $u_-$) divides $\bp_a$ iff $F$ has a root in some torus orbit 
$O_w$ with $w\!\cdot\!a \!>\!0$ (resp.\ $w\cdot a\!<\!0$). More to the 
point, if $c\!\in\!\Ks$ then it follows similarly that $u_++cu_-$ divides 
$\bp_a$ iff some root of $F$ lies in the closure of the hypersurface 
$\{x^a\!=\!c\}$ in $\cT_{P_E}$. Furthermore, lemma 1 immediately 
implies that any two (distinct) hypersurfaces of this form must intersect 
precisely on the ambiguity locus. Thus $\bp_a$ 
can not be identically $0$ and we obtain that $\cN$ is finite. 

Now assume that $F$ has a root in the ambiguity locus. Then, by 
the Vanishing Theorem for Resultants and 
what we've just learned about the hypersurface $\{x^a\!=\!c\}$, 
$\bp_a$ must be divisible by infinitely many distinct binomials of the form 
$u_++cu_-$. So the existence of a root of $F$ within 
the ambiguity locus implies that $\cN\!=\!\infty$. Similarly, 
the existence of infinitely many roots of $F$ within 
$\Ksn$ implies that $\cN\!=\!\infty$ as well.  
Thus we are done with the first portion of the theorem. 
	
As for the codimension of the space of 
``degenerate'' $F$ being $\geq\!2$, first consider the $F$ which have a 
root on the ambiguity locus. These $F$ actually form a codimension $2$ 
subvariety. To see this, first note that $F$ has a root in $O_w 
\Longleftrightarrow$ the initial term system $\init_w(F)$ has a root 
in $\Ksn$ \cite[corollary 2]{toricint}. By the definition of $L(E,a)$, 
we can then conclude via theorem 1.3 of \cite{combiresult}. (Although the 
results of \cite{combiresult} are stated over the complex numbers, this is only 
a minor technicality: in this case, the results of \cite{toricint} immediately 
imply that we can apply theorem 1.3 over {\em any}\/ algebraically closed 
field.) 

As for the $F$ with infinitely many roots in $\Ksn$, assume temporarily that 
the coefficients of $F$ are algebraically independent indeterminates. That 
the $F$ with infinitely many roots in $\Ksn$ define a subvariety of 
codimension $\geq\!2$ must certainly be an old result. However, for 
completeness, we will give a quick proof: 
Following the notation of lemma \ref{lemma:forbid}, Hilbert's 
Irreducibility Theorem \cite{lang} readily implies that 
we can choose $\alpha,\beta,\gamma,\delta\!\in\!\Ks$ such that the 
polynomials $p(\alpha,\beta)$ and $p(\gamma,\delta)$ are relatively 
prime in $K[\cC_E]$. So, by lemma \ref{lemma:forbid}, these $F$ 
are indeed contained in a codimension $2$ subvariety: $\{ 
\cC_E \; | \; p(\alpha,\beta)\!=\!q(\gamma,\delta)\!=\!0 \}.$ \qed 

So in summary, computing $\bp_a$ amounts to laminating $\cT_{P_E}$ with 
hypersurfaces of the form $\overline{\{c_+x^a\!=\!c_-\}}$. These laminae  
all intersect at the ambiguity locus, and toric infinity is precisely the 
union of the two degenerate laminae corresponding to 
(cf.\ remark \ref{rem:egg}) the $\pm$ halves of $\partial P_E$. 

\noindent
{\bf Proof of Theorem 1 and Corollary 1:} We need only state the 
preceding paragraph more algebraically. In particular, lemma \ref{lemma:chow} 
has already done this for us. We thus obtain that when $\cN\!<\!\infty$, 
the multiplicity of a factor $c_+u_++c_-u_-$ of $\bp_a$ is precisely the sum 
of the intersection numbers of the components of $\cD_{P_E}(F,\cP_E)$ lying in 
the fiber $\bar{\varphi}^{-1}([c_-:c_+])$. 
So $\cN$ (resp.\ $\cN'$) is indeed the exact number of roots, counting 
multiplicities (resp.\ {\em not}\/ counting multiplicities), of $F$ 
in $\Ksn$. \qed 

Note that it is actually possible for a positive-dimensional 
component of $\cD_{P_E}(F,\cP_E)$ lying at toric infinity to make a 
positive (finite) contribution to $\cN$ (or $\cN'$). This can be interpreted 
as an alternative way of associating an 
intersection number to an excess component \cite{ifulton}. 



\section{Toric Generalized Characteristic Polynomials} 
\label{sec:toricgcp} 
We will first comment briefly on the combinatorial preprocessing step. The 
following example illustrates a simple fundamental case.
\setcounter{ex}{2}
\begin{ex}[The ``Dense'' Case] 
Suppose $\cP$ is the $n$-tuple $(d_1\Delta,\ldots,d_n\Delta)$ 
where $\Delta\!\subset\!\Rn$ is the standard $n$-simplex and 
$d_i\!\in\!\N$ for all $i$. It is then easily verified that the 
$n$-tuple $D\!:=\!(\{\bO,d_1\hat{e}_1\},\ldots,\{\bO,d_n\hat{e}_n\})$ 
is an irreducible fill of $\cP$ \cite{convexapp,schneider}. Letting $A$ be the 
vertices of $\Delta$, theorem 3 then implies that $\cH$ is a variant (over a 
general algebraically closed field) of the original 
GCP applied to an $n\!\times\!n$ system of equations with degrees 
$d_1,\ldots,d_n$ \cite{gcp}. In particular, our $F-sF^*$ has $2n$ 
$s$-monomials, compared to $n$ $s$-monomials in Canny's 
$(f_1-sx^{d_1}_1,\ldots,f_n-sx^{d_n}_n)$. Note 
also that $\conv(A)$ and $P$ are homothetic and $\cT_P\!\cong\!\Pro^n_K$. 
Neglecting the extra $s$-monomials, setting $d_i\!=\!1$ for all $i$, and 
suitably specializing the coefficients of $g$, we can then recover the usual 
characteristic polynomial of a matrix.
\end{ex} 
We point out that the computational complexity of 
finding an irreducible fill is an open question. However, the 
connection between fills and polynomial system solving (not to mention 
specialized resultants) appears 
to be new and, we hope, provides added incentive to investigate  
filling. Also, even if finding a fill is difficult, this step need only 
be done {\em once}\/ for a given family of problems, provided 
$E$ remains fixed. The situation where the monomial term structure of a 
polynomial system is fixed once and for all (and the coefficients may vary 
thousands of times) actually occurs frequently in many practical contexts, 
such as robot control or computational geometry.

In any event, the toric GCP is an algebraic perturbation method and 
irreducible fills provide a combinatorial means of inducing ``general 
position'' into the roots of $\cF_A$. For example, the following lemma implies 
that $F^*$ is sufficiently generic in a useful sense. 
\begin{lemma}
Following the notation and assumptions of theorem 3, for any point $v$ lying in 
any $D_i$, there exists a $w\!\in\!\Rns$ such that $\{i\}$ is the unique 
essential subset of $D^w$ and $D^w_i\!=\!\{v\}$. In particular, $F^*$ has 
exactly $\cM(\cP)$ roots (counting multiplicities) in $\Ksn$. \qed 
\end{lemma}
\noindent
This lemma follows easily from the techniques of \cite{convexapp}, 
particularly section 2.5. 

Given that $F^*$ thus has maximally many isolated roots (counting 
multiplicities) and {\em no}\/ excess components in $\Ksn$, theorem 3 
then follows easily from the Vanishing Theorem for Resultants and an 
algebraic homotopy argument. 
A similar homotopy technique appears in \cite{rojaswang,
toricint}, and the dense case (with $K\!=\!\C$) is covered in   
\cite{gcp,mikebezout}. (However, toric variety language was not used in the 
last two works.) Thus, the usual GCP is (almost) 
the special case of the toric GCP where $\cT_P$ is complex projective space. 

In closing, an important difference to note is that our present toric GCP is 
primarily suited for $\Ksn$, while the original GCP is mainly suited (in a 
non-sparse way) for affine space. To {\em completely}\/ generalize and improve 
the GCP in affine space, it is necessary 
to use the {\em affine}\/  sparse resultant and this is pursued further in 
\cite{aff}. For instance, by replacing the sparse resultant with the affine 
sparse resultant, and using {\em $\Kn$-counting}\/ \cite{toricint} instead of 
filling, we can actually recover Canny's GCP in the dense case. 

\section{Acknowledgements} 
The main results of this paper were presented 
during a visit of the author at City University of Hong Kong. The author 
thanks Steve Smale for his invitation, warm hospitality, and even warmer beach 
hikes. Special thanks go to John Dalbec for pointing on an error in an 
earlier version of section 3, and to Robert L. Williams for his encouragement 
and keen commentary. The author also thanks his former roommates 
Jean-Pierre Dedieu and Gregorio Malajovich for their input, but most 
of all for not getting too mad at him for not washing the dishes. 

Work on this paper began at MSRI and concluded at MIT. The author 
thanks both institutions for their support, and especially Marsh 
Borg at MSRI for her kind assistance. 
 
\bibliographystyle{amsalpha}

\end{document}